\documentclass[12pt]{article} \usepackage{amsfonts}
\begin{document}
\baselineskip=14pt
\pagestyle{plain}
{\Large

\vspace{10mm}

\centerline {\bf On the basis property of the root function systems of regular}
\centerline {\bf boundary value problems for the Sturm-Liouville operator}

\medskip
\medskip
\centerline { Alexander Makin}
\vspace{10mm}

\centerline{Abstract. We consider the nonselfadjoint Sturm-Liouville}
\centerline{operator with regular but not strongly regular boundary}
\centerline{conditions. We examine the basis property of the root}
\centerline{function system of the mentioned operator.}
\vspace{10mm}

In the present paper we study eigenvalue problems for the nonselfadjoint
Sturm-Liouville operator
$$
Lu=u''-q(x)u \eqno (1)
$$
defined on the interval $(0,1)$, where $q(x)$ is an arbitrary complex-valued
function of the class $L_1(0,1)$. Our main purpose is to investigate the basis
property of the root function system of operator (1) with regular but not strongly
regular boundary conditions. Author's interest to this problem was stimulated
by the papers of V.A. Il'in [1-3]. 

By $\varphi(x,\mu), \psi(x,\mu)$ we denote the fundamental for $\mu\ne0$
system of solutions to the equation
$$
u''-q(x)u+\mu^2u=0 $$
determined by the initial conditions $\varphi(0,\mu)=\psi(0,\mu)=1$, $\varphi_x'(0,\mu)=i\mu$,
$\psi_x'(0,\mu)=-i\mu$. It is well known that the functions $\varphi(x,\mu)$ and
$\psi(x,\mu)$ satisfy the integral equations
$$
\varphi(x,\mu)=e^{i\mu x}+\frac{1}{\mu}\int_0^x\sin\mu(x-t)q(t)\varphi(t,\mu)dt, \eqno (2)
$$
$$
\psi(x,\mu)=e^{-i\mu x}+ \frac{1}{\mu}\int_0^x\sin\mu(x-t)q(t)\psi(t,\mu)dt, \eqno(2')
$$
respectively. Also it is well known that these functions are continuous with
their partial derivatives, and for any fixed $x$ they are analytic functions
of the parameter $\mu$. Later on, we suppose that the inequality

$$
|Im\mu|<M \eqno(3)
$$
holds, where $M$ is some constant. It is known [4] that the estimates
$$
|\varphi(x,\mu)|\le c_1, \qquad |\varphi_\mu'(x,\mu)|\le c_1, \eqno(4)
$$
$$
|\psi(x,\mu)|\le c_2, \qquad |\psi_\mu'(x, \mu)|\le c_2 \eqno(4')
$$
are valid for $0\le x\le1$ and $|\mu|\ge\mu_0$, where $\mu_0$ is a
sufficiently large number.

We need more precise  asymptotic formulas for the functions $\varphi(x,\mu)$
and $\psi(x,\mu)$. Transforming the right-hand side of (2), we have
$$
\begin{array}{c}
\varphi(x,\mu)=e^{i\mu x}+(2i\mu)^{-1}\int_0^x(e^{i\mu(x-t)}-e^{-i\mu(x-t)})q(t)\times\\
\times[e^{i\mu t}+\mu^{-1}\int_0^t\sin\mu(t-s)q(s)\varphi(t,\mu)ds]dt=\\
=e^{i\mu x}\{1+(2i\mu)^{-1}\int_0^xq(t)dt-(2i\mu)^{-1}\int_0^xe^{2i\mu(t-x)}q(t)dt-\\
-(4\mu^2)^{-1}\int_0^x(e^{i\mu(x-t)}-e^{-i\mu(x-t)})q(t)dt\times\\
\times\int_0^t(e^{i\mu(t-s)}-e^{-i\mu(t-s)})q(s)[e^{i\mu s}+\mu^{-1}\int_0^s\sin\mu(s-y)q(y)\varphi(y,\mu)dy]ds\}.
\end{array}\eqno(5)
$$
We denote the sum of the first three summands in braces on the right-hand side
of (5) by $F_1(x,\mu)$. Dividing the integral with respect to $s$ on the right-hand side of (5) into three
summands, we obtain
$$
\begin{array}{c}
\varphi(x,\mu)=e^{i\mu x}\{F_1(x,\mu)-(4\mu^2)^{-1}e^{-i\mu x}\int_0^x(e^{i\mu(x-t)}-e^{-i\mu(x-t)})q(t)dt\times\\
\times[e^{i\mu t}\int_0^tq(s)ds-\int_0^te^{-i\mu(t-2s)}q(s)ds+\\
+\mu^{-1}\int_0^t(e^{i\mu(t-s)}-e^{-i\mu(t-s)})q(s)ds\int_0^s\sin\mu(s-y)q(y)\varphi(y,\mu)dy]\}.
\end{array}\eqno(6)
$$
We denote the last summand in square brackets on the right-hand side of (6) by $\theta_1(t,\mu)$.
It follows from (3) and (4) that
$
\theta_1(t,\mu)=O(\mu^{-1})$ and $ \partial\theta_1(t,\mu)/\partial\mu=O(\mu^{-1})$.

Writing the terms in braces on the right-hand side of (6) in descending powers of $\mu$,
we get
$$
\begin{array}{c}
\varphi(x,\mu)=e^{i\mu x}\{F_1(x,\mu)-(4\mu^2)^{-1}\int_0^xq(t)dt\int_0^tq(s)ds+\\
+(4\mu^2)^{-1}\int_0^xe^{2i\mu(t-x)}q(t)dt\int_0^tq(s)ds+\\
+(4\mu^2)^{-1}e^{-i\mu x}\int_0^x(e^{i\mu(x-t)}-e^{-i\mu(x-t)})q(t)dt\int_0^te^{-i\mu(t-2s)}q(s)ds-\\
-(4\mu^2)^{-1}e^{-i\mu x}\int_0^x(e^{i\mu(x-t)}-e^{-i\mu(x-t)})q(t)\theta_1(t,\mu)dt\}.
\end{array}\eqno(7)
$$
We denote the last summand in braces on the right-hand side of (7) by $\theta_2(x,\mu)$.
It is readily seen that
$
\theta_2(x,\mu)=O(\mu^{-3})$ and $ \partial\theta_2(x,\mu)/\partial\mu=O(\mu^{-3})$.
Simplifying the expression in braces on the right-hand side of (7), we obtain
$$
\begin{array}{c}
\varphi(x,\mu)=e^{i\mu x}\{F_1(x,\mu)-(4\mu^2)^{-1}\int_0^xq(t)dt\int_0^tq(s)ds+\\
+[(4\mu^2)^{-1}(e^{-2i\mu x}\int_0^xe^{2i\mu t}q(t)dt\int_0^tq(s)ds+\int_0^xe^{-2i\mu t}q(t)dt\int_0^te^{2i\mu s}q(s)ds-\\
-e^{-2i\mu x}\int_0^xq(t)dt\int_0^te^{2i\mu s}q(s)ds)]+\theta_2(x,\mu)\}.
\end{array}\eqno(8)
$$ 
We denote the expression in square brackets on the right-hand side of (8) by $\theta_3(x,\mu)$.
It follows from (3) and the Riemann lemma that
$
\theta_3(x,\mu)=o(\mu^{-2})$ and $ \partial\theta_3(x,\mu)/\partial\mu=o(\mu^{-2})$.

Transforming the right-hand side of (2$'$), we have
$$
\begin{array}{c}
\psi(x,\mu)=e^{-i\mu x}+(2i\mu)^{-1}\int_0^x(e^{i\mu(x-t)}-e^{-i\mu(x-t)})q(t)\times\\
\times[e^{-i\mu t}+\mu^{-1}\int_0^t\sin\mu(t-s)q(s)\psi(s,\mu)ds]dt=\\
=e^{-i\mu x}\{1-(2i\mu)^{-1}\int_0^xq(t)dt+(2i\mu)^{-1}\int_0^xe^{2i\mu(x-t)}q(t)dt-\\
-(4\mu^2)^{-1}e^{i\mu x}\int_0^x(e^{i\mu(x-t)}-e^{-i\mu(x-t)})q(t)dt\times\\
\times\int_0^t(e^{i\mu(t-s)}-e^{-i\mu(t-s)})q(s)[e^{-i\mu s}+\mu^{-1}\int_0^s\sin\mu(s-y)q(y)\psi(y,\mu)dy]ds\}.
\end{array}\eqno(5')
$$
We denote the sum of the first three summands in braces on the right-hand side of (5$'$)
by $F_2(x,\mu)$. Dividing the integral with respect to $s$ on the right-hand side of (5$'$) into the sum of three
summands, we obtain
$$
\begin{array}{c}
\psi(x,\mu)=e^{-i\mu x}\{F_2(x,\mu)-(4\mu^2)^{-1}e^{i\mu x}\int_0^x(e^{i\mu(x-t)}-e^{-i\mu(x-t)})q(t)dt\times\\
\times[-e^{-i\mu t}\int_0^tq(t)dt+\int_0^te^{i\mu(t-2s)}q(s)ds+\\
+\mu^{-1}\int_0^t(e^{i\mu(t-s)}-e^{-i\mu(t-s)})q(s)ds\int_0^s\sin\mu(s-y)q(y)\psi(y,\mu)dy]\}.
\end{array}\eqno(6')
$$
We denote the last summand in square brackets on the right-hand side of (6$'$) by $\tilde\theta_1(t,\mu)$.
It follows from (3) and (4$'$) that
$
\tilde\theta_1(t,\mu)=O(\mu^{-1})$ and $ \partial\tilde\theta_1(t,\mu)/\partial\mu=O(\mu^{-1})$.

Writing the terms in braces on the right-hand side of (6$'$) in decsending powers of $\mu$,
we get
$$
\begin{array}{c}
\psi(x,\mu)=e^{-i\mu x}\{F_2(x,\mu)-(4\mu^2)^{-1}\int_0^xq(t)dt\int_0^tq(s)ds+\\
+(4\mu^2)^{-1}\int_0^xe^{2i\mu(x-t)}q(t)dt\int_0^tq(s)ds-\\
-(4\mu^2)^{-1}e^{i\mu x}\int_0^x(e^{i\mu(x-t)}-e^{-i\mu(x-t)})q(t)dt\int_0^te^{i\mu(t-2s)}q(s)ds-\\
-(4\mu^2)^{-1}e^{i\mu x}\int_0^x(e^{i\mu(x-t)}-e^{-i\mu(x-t)})q(t)\tilde\theta_1(t,\mu)dt\}.
\end{array}\eqno(7')
$$
We denote the last summand in braces on the right-hand side (7$'$) by $\tilde\theta(x,\mu)$.
It is readily seen that 
$
\tilde\theta_2(x,\mu)=O(\mu^{-3})$ and $ \partial\tilde\theta_2(x,\mu)/\partial\mu=O(\mu^{-3})$.

Simplifying the expression in braces on the right-hand side of (7$'$), we obtain
$$
\begin{array}{c}
\psi(x,\mu)=e^{-i\mu x}\{F_2(x,\mu)-(4\mu^2)^{-1}\int_0^xq(t)dt\int_0^tq(s)ds+\\
+[(4\mu^2)^{-1}(e^{2i\mu x}\int_0^xe^{-2i\mu t}dt\int_0^tq(s)ds+\int_0^xe^{2i\mu t}q(t)dt\int_0^te^{-2i\mu s}q(s)ds-\\
-e^{2i\mu x}\int_0^xq(t)dt\int_0^te^{-2i\mu s}q(s)ds)]+\tilde\theta_2(x,\mu)\}.
\end{array}\eqno(8')
$$
We denote the expression in square brackets on the right-hand side (8$'$) by $\tilde\theta_3(x,\mu)$.
It follows from (3) and the Riemann lemma that
$
\tilde\theta_3(x,\mu)=o(\mu^{-2})$ and $ \partial\tilde\theta_3(x,\mu)/\partial\mu=o(\mu^{-2})$.

In the same way, we get asymptotic formulas for the functions $\varphi_x'(x,\mu)$
and $\psi_x'(x,\mu)$. Differentiating relations (2) and (2$'$), we have
$$
\varphi_x'(x,\mu)=i\mu e^{i\mu x}+\int_0^x\cos\mu(x-t)q(t)\varphi(t,\mu)dt,\eqno(9)
$$
$$
\psi_x'(x,\mu)=-i\mu e^{-i\mu}+\int_0^x\cos\mu(x-t)q(t)\psi(t,\mu)dt.\eqno(9')
$$
Transforming the second summand on the right-hand side of (9) according to (8), we have
$$
\begin{array}{c}
\varphi_x'(x,\mu)=i\mu e^{i\mu x}+\frac{1}{2}\int_0^x(e^{i\mu(x-t)}+e^{-i\mu(x-t)})q(t)\times\\
\times[e^{i\mu t}(F_1(t,\mu)+\theta_4(t,\mu)]dt=\\
=i\mu e^{i\mu x}+\frac{1}{2}\int_0^x(e^{i\mu x}+e^{i\mu(2t-x)})q(t)
[F_1(t,\mu)+
\theta_4(t,\mu)]dt,
\end{array}\eqno(10)
$$
where
$$
\theta_4(t,\mu)=-(4\mu^2)^{-1}\int_0^tq(s)ds\int_0^sq(y)dy+\theta_2(t,\mu)+\theta_3(t,\mu).
$$
Evidently,
$
\theta_4(t,\mu)=O(\mu^{-2})$ and $ \partial\theta_4(t,\mu)/\partial\mu=O(\mu^{-2})$.

Writing the terms on the right-hand side of (10) in decsending powers of $\mu$,
we get
$$
\begin{array}{c}
\varphi_x'(x,\mu)=e^{i\mu x}\{i\mu+\frac{1}{2}\int_0^xq(t)dt+\frac{1}{2}\int_0^xe^{2i\mu(t-x)}q(t)dt+\\
+(4i\mu)^{-1}\int_0^xq(t)dt\int_0^tq(s)ds+\\
+[-(4i\mu)^{-1}\int_0^xe^{-2i\mu t}q(t)dt\int_0^te^{2i\mu s}q(s)ds+\\
+(4i\mu)^{-1}e^{-2i\mu x}\int_0^xe^{2i\mu t}q(t)dt\int_0^tq(s)ds-\\
-(4i\mu)^{-1}e^{-2i\mu x}\int_0^xq(t)dt\int_0^te^{2i\mu s}q(s)ds]+\\
+\frac{1}{2}\int_0^x(1+e^{2i\mu(t-x)})q(t)\theta_4(t,\mu)dt\}.
\end{array}\eqno(11)
$$
We denote the expression in square brackets of (11) by $\theta_5(x,\mu)$. It follows from (3) and
the Riemann lemma that
$
\theta_5(x,\mu)=o(\mu^{-1})$ and $ \partial\theta_5(x,\mu)/\partial\mu=o(\mu^{-1})$.

Trasforming the second summand on the right-hand side of (9$'$) according to (8$'$), we obtain
$$
\begin{array}{c}
\psi_x'(x,\mu)=-i\mu e^{-i\mu x}+\frac{1}{2}\int_0^x(e^{i\mu(x-t)}+e^{-i\mu(x-t)})q(t)\times\\
\times[e^{-i\mu t}(F_2(t,\mu)+\tilde\theta(t,\mu)]dt=\\
=-i\mu e^{-i\mu x}+\frac{1}{2}\int_0^x(e^{i\mu(x-2t)}+e^{-i\mu x})q(t)[F_2(t,\mu)+\tilde\theta_4(t,\mu)]dt,
\end{array}\eqno(10')
$$
where
$$
\tilde\theta_4(t,\mu)=-(4\mu^2)^{-1}\int_0^tq(s)ds\int_0^sq(y)dy+\tilde\theta_2(t,\mu)+\tilde\theta_3(t,\mu).
$$
Obviously,
$$
\tilde\theta_4(t,\mu)=O(\mu^{-2}), \qquad \partial\tilde\theta_4(t,\mu)/\partial\mu=O(\mu^{-2}).
$$
Writing the terms on the right-hand side of (10$'$) in decsending powers of $\mu$,
we get
$$
\begin{array}{c}
\psi_x'(x,\mu)=e^{-i\mu x}\{-i\mu+\frac{1}{2}\int_0^xq(t)dt+\frac{1}{2}\int_0^xe^{2i\mu(x-t)}q(t)dt-\\
-(4i\mu)^{-1}\int_0^xq(t)dt\int_0^tq(s)ds+\\
+[(4i\mu)^{-1}\int_0^xe^{2i\mu t}q(t)dt\int_0^te^{-2i\mu s}q(s)ds-\\
-(4i\mu)^{-1}e^{2i\mu x}\int_0^xe^{-2i\mu t}q(t)dt\int_0^tq(s)ds+\\
+(4i\mu)^{-1}e^{2i\mu x}\int_0^xq(t)dt\int_0^te^{-2i\mu s}q(s)ds]+\\
+\frac{1}{2}\int_0^x(e^{2i\mu(x-t)}+1)q(t)\tilde\theta_4(t,\mu)\}.
\end{array}\eqno(11')
$$
We denote the expression in square brackets in (11$'$) by $\tilde\theta_5(x,\mu)$. It folows from (3) and the
Riemann lemma that
$$
\tilde\theta_5(x,\mu)=o(\mu^{-1}), \qquad \partial\tilde\theta_5(x,\mu)/\partial\mu=o(\mu^{-1}).
$$

To simplify the asymptotic formulas obtained above, let us prove that
$$
\int_0^xq(t)dt\int_0^tq(s)ds=\frac{1}{2}\left(\int_0^xq(t)dt\right)^2.\eqno(12)
$$
If the function $q(t)$ is continuous on the segment $[0,1]$, then equality (12)
can be obtained by integration by parts. In general case, let us
approximate the function $q(t)$ by a continuous function $f(t)$ so that $q(t)-f(t)=r(t)$,
where $\int_0^1|r(t)|dt<\varepsilon$, where $\varepsilon>0$ is an arbitrary preassigned
number. Denote $q_0=\int_0^1|q(t)|dt$. Then we have 
$$
\begin{array}{c}
|\int_0^xq(t)dt\int_0^tq(s)ds-\frac{1}{2}\left(\int_0^xq(t)dt\right)^2|=\\
=|\int_0^x(f(t)+r(t))dt\int_0^t(f(s)+r(s))ds-\frac{1}{2}\left(\int_0^x(f(t)+r(t))dt\right)^2|=\\
=|\int_0^xf(t)dt\int_0^tf(s)ds+\int_0^xr(t)dt\int_0^tf(s)ds+\int_0^xf(t)dt\int_0^tr(s)ds+\\
+\int_0^xr(t)dt\int_0^tr(s)ds-\frac{1}{2}\left(\int_0^xf(t)dt\right)^2-\int_0^xf(t)dt\int_0^xr(t)dt-\\
-\frac{1}{2}\left(\int_0^xr(t)dt\right)^2|\le\varepsilon(3q_0+3\varepsilon/2).
\end{array}
$$
This yields that equality (12) holds for any function $q(t)$ of the class
$L_1(0,1)$. It follows from (8), (8$'$), (11), (11$'$), the estimates for the functions $\theta_i$,
$\tilde\theta_i$ $(i=\overline{1,5})$ and (12) that
$$
\begin{array}{c}
\varphi(x,\mu)=e^{i\mu x}\{1+(2i\mu)^{-1}\int_0^xq(t)dt-(2i\mu)^{-1}\int_0^xe^{2i\mu(t-x)}q(t)dt-\\
-(8\mu^2)^{-1}\left(\int_0^xq(t)dt\right)^2+\theta_6(x,\mu)\},
\end{array}\eqno(13)
$$
$$
\begin{array}{c}
\psi(x,\mu)=e^{-i\mu x}\{1-(2i\mu)^{-1}\int_0^xq(t)dt+(2i\mu)^{-1}\int_0^xe^{2i\mu(x-t)}q(t)dt-\\
-(8\mu^2)^{-1}\left(\int_0^xq(t)dt\right)^2+\tilde\theta_6(x,\mu)\},
\end{array}\eqno(13')
$$
where
$
\theta_6(x,\mu)=o(\mu^{-2})$, $\partial\theta_6(x,\mu)/\partial\mu=o(\mu^{-2})$,
$\tilde\theta_6(x,\mu)=o(\mu^{-2})$, $\partial\tilde\theta_6(x,\mu)/\partial\mu=o(\mu^{-2})$;

$$
\begin{array}{c}
\varphi_x'(x,\mu)=e^{i\mu x}\{i\mu+\frac{1}{2}\int_0^xq(t)dt+\frac{1}{2}\int_0^xe^{2i\mu(t-x)}q(t)dt+\\
+(8i\mu)^{-1}\left(\int_0^xq(t)dt\right)^2+\theta_7(x,\mu)\},
\end{array}\eqno(14)
$$
$$
\begin{array}{c}
\psi_x'(x,\mu)=e^{-i\mu x}\{-i\mu+\frac{1}{2}\int_0^xq(t)dt+\frac{1}{2}\int_0^xe^{2i\mu(x-t)}q(t)dt-\\
-(8i\mu)^{-1}\left(\int_0^xq(t)dt\right)^2+\tilde\theta_7(x,\mu)\},
\end{array}\eqno(14')
$$
where
$
\theta_7(x,\mu)=o(\mu^{-1})$, $\partial\theta_7(x,\mu)/\partial\mu=o(\mu^{-1})$,
$\tilde\theta_7(x,\mu)=o(\mu^{-1})$, $\partial\tilde\theta_7(x,\mu)/\partial\mu=o(\mu^{-1})$.

For operator (1) let us consider the following two-point
boundary value problem with boundary conditions determined by linearly independent
forms with arbitrary complex-valued coefficients
$$
\begin{array}{c}
B_1(u)=a_1u'(0)+b_1u'(1)+a_0u(0)+b_0u(1)=0,\\ 
B_2(u)=c_1u'(0)+d_1u'(1)+c_0u(0)+d_0u(1)=0.
\end{array}\eqno(15)
$$
It is convenient to rewrite conditions (15) in terms of the matrix $A$, where
$$
A=\left(
\begin{array}{cccc}
a_1&b_1&a_0&b_0\\
c_1&d_1&c_0&d_0
\end{array}
\right);
$$
by $A(ij)$ we denote the matrix consisting of the $i$th and $j$th
columns of the matrix $A$ $(1\le i\le j\le4)$, and we set $A_{ij}=det$ $A(ij)$. We also denote
$$
B_1^\star(u)=-a_1u'(0)+b_1u'(1)-a_0u(0)+b_0u(1),
$$
$$
B_2^\star(u)=-c_1u'(0)+d_1u'(1)-c_0u(0)+d_0u(0).
$$

In addition, we assume that
$$
\int_0^1q(x)dx=0. \eqno(16)
$$
Then it follows from (13), (13$'$), (14), (14$'$) and (16) that
$$
\begin{array}{c}
\varphi(0,\mu)=1, \qquad\varphi(1,\mu)=e^{i\mu}(1+P),\\
\varphi_x'(0,\mu)=i\mu, \qquad\varphi_x'(1,\mu)=e^{i\mu}(i\mu+P'),
\end{array}\eqno(17)
$$
where
$$
\begin{array}{c}
P=-(2i\mu)^{-1}e^{-2i\mu}\int_0^1e^{2i\mu t}q(t)dt+\theta_6(1,\mu),\\
P'=\frac{1}{2}e^{-2i\mu}\int_0^1e^{2i\mu t}q(t)dt+\theta_7(1,\mu);
\end{array}\eqno(18)
$$
$$
\begin{array}{c}
\psi(0,\mu)=1, \qquad\psi(1,\mu)=e^{-i\mu}(1+Q),\\
\psi_x'(0,\mu)=-i\mu, \qquad\psi_x'(1,\mu)=e^{-i\mu}(-i\mu+Q'),
\end{array}\eqno(17')
$$
where
$$
\begin{array}{c}
Q=(2i\mu)^{-1}e^{2i\mu}\int_0^1e^{-2i\mu t}q(t)dt+\tilde\theta_6(1,\mu),\\
Q'=\frac{1}{2}e^{2i\mu}\int_0^1e^{-2i\mu t}q(t)dt+\tilde\theta_7(1,\mu).
\end{array}\eqno(18')
$$

By $\Delta(\mu)$ we denote the characteristic determinant of the problem
$$
Lu+\mu^2u=0, \quad B_1(u)=0, \quad B_2(u)=0,\eqno(19)
$$
and by $\Delta_0(\mu)$ we denote the characteristic determinant of the problem
$$
u''+\mu^2u=0, \quad B_1(u)=0, \quad B_2(u)=0.
$$
Using relations (17-18$'$), and performing some simple though awkward 
manipulations, we obtain
$$
\begin{array}{c}
\Delta(\mu)=B_1(\varphi)B_2(\psi)-B_1(\psi)B_2(\varphi)=\\
=[a_1i\mu+b_1e^{i\mu}(i\mu+P')+a_0+b_0e^{i\mu}(1+P)]\times\\
\times[c_1(-i\mu)+d_1e^{-i\mu}(-i\mu+Q')+c_0+d_0e^{-i\mu}(1+Q)]-\\
-[a_1(-i\mu)+b_1e^{-i\mu}(-i\mu)(-i\mu+Q')+a_0+b_0e^{-i\mu}(1+Q)]\times\\
\times[c_1i\mu+d_1e^{i\mu}(i\mu+P')+c_0+d_0e^{i\mu}(1+P)]=\\
=\Delta_0(\mu)+i\mu A_{12}(e^{i\mu}P'+e^{-i\mu}Q')+i\mu A^{14}(e^{i\mu}P+e^{-i\mu}Q)+\\
+A_{23}(e^{i\mu}P'-e^{-i\mu}q')+A_{24}[P'-Q'+i\mu(P+Q)]+\alpha_0(\mu),
\end{array}
$$
where
$$
\begin{array}{c}
\alpha_0(\mu)= b_1d_1P'Q'+b_1d_0P'Q+b_0d_1PQ'+b_0d_0PQ-\\
-(b_1d_1Q'P'+b_1d_0Q'P+b_0d_1QP'+b_0d_0QP)=\\
=b_1d_0(P'Q-PQ')+b_0d_1(PQ'-P'Q)=A_{24}(P'Q-PQ').
\end{array}
$$
It follows from (18) and (18$'$) that $\alpha_0(\mu)=o(\mu^{-1})$ and $\alpha_0'(\mu)=o(\mu^{-1})$.
In the following, we assume that $A_{12}=0$. It follows from (18) and (18$'$) that
$$
\begin{array}{c}
i\mu A_{14}(e^{i\mu}P+e^{-i\mu}Q)=\\
=\frac{1}{2}A_{14}(e^{i\mu}\int_0^1e^{-2i\mu t}q(t)dt-e^{-i\mu}\int_0^1e^{2i\mu t}q(t)dt)+\alpha_1(\mu),\\
A_{23}(e^{i\mu}P'-e^{-i\mu}Q')=\\
=\frac{1}{2}A_{23}(e^{-i\mu}\int_0^1e^{2i\mu t}q(t)dt-e^{i\mu}\int_0^1e^{-2i\mu t}q(t)dt)+\alpha_2(\mu),\\
A_{24}(P'-Q'+i\mu(P+Q))=\alpha_3(\mu),
\end{array}
$$
where $\alpha_j(\mu)=o(\mu^{-1})$ and $\alpha_j'(\mu)=o(\mu^{-1})$, $j=1,2,3$. Hence, 
$$
\Delta(\mu)=\Delta_0(\mu)+\frac{1}{2}(A_{14}-A_{23})[e^{i\mu}\int_0^1e^{-2i\mu t}q(t)dt-e^{-i\mu}\int_0^1e^{2i\mu t}q(t)dt]+\theta(\mu), \eqno(20)
$$
where $\theta(\mu)=o(\mu^{-1})$ and $\theta'(\mu)=o(\mu^{-1})$.

Let boundary conditions (15) be regular but not strongly regular [5, pp. 71-73],
which, by [5, p. 73] is equivalent to the conditions
$$
A_{12}=0, \quad A_{14}+A_{23}\ne0,\quad A_{14}+A_{23}=\mp(A_{13}+A_{24}).\eqno(21)
$$
Without loss of generality, we assume that conditions (15) are normalized [5, p. 66]. This, together
with the relation $A_{12}=0$, yields
$$
c_1=d_1=0.\eqno(22)
$$

Let $\{u_n(x)\}$ be the system of eigenfunctions and associated functions of problem (19),
and let $\lambda_n=\mu_n^2$ be the corresponding eigenvalues $(Re\mu_n\ge0)$. By [5, p. 74] the set of numbers $\mu_n$, except for
possibly finitely many numbers, consists of two series
$$
\mu_n'=2\pi n+\delta_n',\quad \mu_n''=2\pi n+\delta_n',\eqno(23)
$$
where $|\delta_n'|\le c_1n^{-1/2}$, $|\delta_n''|\le c_1n^{-1/2}$, $c_1>0$, $n=n_0, n_0+1, \ldots$,
if $A_{14}+A_{23}=-(A_{13}+A_{24})$ in (21) (case 1), and this set, except
for possibly finitely many numbers, consists of two series
$$
\tilde\mu_n'=(2n-1)\pi+\tilde\delta_n',\quad \tilde\mu_n''=(2n-1)\pi+\tilde\delta_n'',\eqno(23')
$$
where $|\tilde\delta_n'|\le c_3n^{-1/2}$, $|\tilde\delta_n''|\le c_3n^{-1/2}$, $c_3>0$,
$n=n_0, n_0+1,\ldots$, if $A_{14}+A_{23}=A_{13}+A_{24}$ in (21) (case 2).
It is also known 
[5, p. 98, p. 91] that the system $\{u_n(x)\}$ is complete in $L_2(0,1)$ and there
exists a biorthogonally conjugate system $\{v_n(x)\}$.
If $\mu_n$ is a simple zero of the function
$\Delta(\mu)$, then, by [5, p. 48],
$$
u_n(x)\overline{v_n(\xi)}=-2\mu_nH(x,\xi,\mu_n)/\Delta'(\mu_n),\eqno(24)
$$
where
$$
H(x,\xi,\mu)=\left|
\begin{array}{ccc}
\varphi(x)&\psi(x)&g(x,\xi)\\
B_1(\varphi)&B_1(\psi)&B_1(g)\\
B_2(\varphi)&B_2(\psi)&B_2(g)
\end{array}
\right|,\eqno(25)
$$
$$
g(x,\xi)=\pm\frac{1}{2W(\xi)}\left|
\begin{array}{cc}
\varphi(x)&\psi(x)\\
\varphi(\xi)&\psi(\xi)
\end{array}
\right|,\quad W(\xi)=\left|
\begin{array}{cc}
\varphi'(\xi)&\psi'(\xi)\\
\varphi(\xi)&\psi(\xi)
\end{array}
\right|;\eqno(26)
$$
moreover, the sign "+" corresponds to the case $x>\xi$, and the sign "-" corresponds to the case $x<\xi$. Developing 
determinants (25) and (26), we obtain
$$
H(x,\xi,\mu_n)=\Phi(x,\xi,\mu_n)/(2W(\xi)),\eqno(27)
$$
where
$$
\begin{array}{c}
\Phi(x,\xi,\mu)=\varphi(x)[B_1(\psi)(\psi(\xi)B_2^{\star}(\varphi)-\varphi(\xi)B_2^{\star}(\psi))-\\
-B_2(\psi)(\psi(\xi)B_1^{\star}(\varphi)-\varphi(\xi)B_1^{\star}(\psi))]-\psi(x)[B_1(\varphi)(\psi(\xi)B_2^{\star}(\varphi)-\\
-\varphi(\xi)B_2^{\star}(\psi))-B_2(\varphi)(\psi(\xi)B_1^{\star}(\varphi)-\varphi(\xi)B_1^{\star}(\psi))].
\end{array}\eqno(28)
$$

{\bf Theorem 1.} {\it If $A_{14}=A_{23}$ and $A_{34}\ne0$, then the system $\{u_n(x)\}$
forms a Riesz basis in $L_2(0,1)$.}

Proof. It follows from [6] that
$$
\Delta_0(\mu)=\mp i(A_{13}+A_{24})\mu e^{-i\mu}(e^{i\mu}\mp1)[(e^{i\mu}\mp1)\pm\frac{A_{34}}{i(A_{13}+A_{24})\mu}(e^{i\mu}\pm1)],
$$
where the upper sign is chosen in the case 1
and the lower sign is chosen in the case 2.
This, together with (20), reduces the equation $\Delta(\mu)=0$ to the form 
$$
\mp i(A_{13}+A_{24})\mu(1\mp e^{-i\mu})[(e^{i\mu}\mp1)\pm\frac{A_{34}}{i(A_{13}+A_{24})\mu}(e^{i\mu}\pm1)]+\theta(\mu)=0.
$$
The last equation is reduced to the form
$$
w_1(\mu)=(1-e^{-i\mu})[(e^{i\mu}-1)+\frac{b}{i\mu}(e^{i\mu}+1)]+R_1(\mu)=0\eqno(29)
$$
in the case 1, and
$$
w_2(\mu)=(1+e^{-i\mu})[(e^{i\mu}+1)-\frac{b}{i\mu}(e^{i\mu}-1)]+R_2(\mu)=0\eqno(29')
$$
in the case 2, where $b=A_{34}/(A_{13}+A_{24})$, $R_j(\mu)=o(\mu^{-2})$, $R_j'(\mu)=o(\mu^{-2})$, $j=1,2$.

Let us consider case 1.
Substituting $\mu=2\pi n+z$ into (29) and using (23) we find that the function
$F_n(z)=g(z)f_n(z)+R_1(2\pi n+z)$, where $g(z)=1-e^{-iz}$, $f_n(z)=e^{iz}-1+b(e^{iz}+1)/(i(2\pi n+z))$,
has two roots $\delta_n'$ and $\delta_n''$ in the disk $|z|\le c_1n^{-1/2}$. Evidently,
the function $g(z)$ has a unique root $z=0$ in the same disk, moreover, it follows from [6]
that the function $f_n(z)$ has a unique root $z_n''$ in the same disk,
and $z_n''=O(n^{-1})$. It follows from the last equality and the Maclaurin formula for the function
$e^{iz}$ that $z_n''=b/(\pi n)+O(n^{-2})$.

By $\Gamma_n'$ and $\Gamma_n''$ we denote the circles of radius $r_n=|b|/(4\pi n)$ centered
at $0$ and $b/(\pi n)$, respectively. It follows from the Maclaurin formula that
for all sufficiently large
 $n$ for $z\in\Gamma_n'\cup \Gamma_n''$ $|g(z)f_n(z)|\ge c_2n^{-2}$ $(c_2>0)$. Therefore,
for all sufficiently large $n$ for $z\in \Gamma_n'\cup \Gamma_n''$
$|g(z)f_n(z)|>|F_n(z)-g(z)f_n(z)|$. By the Rouche' theorem, it follows from the last
inequality that the functions $g(z)f_n(z)$ and $F_n(z)$ have the same number of zeros 
in the disks bounded by $\Gamma_n'$ and $\Gamma_n''$, hence, for all sufficiently large $n$ the equation $F_n(z)=0$ has exactly one root
in each disk bounded by $\Gamma_n'$ or $\Gamma_n''$. Thus, we have 
$$
|\delta_n'|<r_n, \quad|\delta_n''-b/(\pi n)|<r_n, \quad|\mu_n'-\mu_n''|>2r_n. \eqno(30)
$$

 In case 2 equation (29$'$) can be analysed in a similar way. Arguing as above, we see that
$$
|\tilde\delta_n'|<r_n, \quad|\tilde\delta_n''+b/(\pi n)|<r_n, \quad|\tilde\mu_n'-\tilde\mu_n''|>2r_n \eqno(30')
$$
for all sufficiently large $n$. In particular, it follows from (30) and (30$'$) that
the eigenvalues $\lambda_n$ are asymptotically simple.

Let us prove that for all sufficiently large $n$ 
$$
c_4\le|\Delta'(\mu_n)|\le c_5,\eqno(31)
$$ 
where $c_4>0$ and $c_5>0$; here $\mu_n$ is an arbitrary root of the equation $\Delta(\mu)=0$.
In case 1 we have $\Delta(\mu)=\beta_1\mu w_1(\mu)$, where $\beta_1=-i(A_{13}+A_{24})$,
therefore,
$$
\Delta'(\mu_n)=\beta_1\mu_nw_1'(\mu_n).\eqno(32)
$$
Let us estimate the function $w_1'(\mu)$. If $\mu=2\pi n+z$, then $w_1(\mu)=F_n(z)$. It
follows from (30) and the Maclaurin formula that for $z=\delta_n'$ we have
$c_6\le|g'(z)|\le c_7$, $c_8/n\le|f_n(z)|\le c_9/n$, $|g(z)|\le c_{10}/n$,
$|f_n(z)|\le c_{11}/n$, and for $z=\delta_n''$ we have $c_{12}\le|f_n(z)|\le c_{13}$,
$c_{14}/n\le|g(z)|\le c_{15}/n$, $|f_n(z)|\le c_{16}/n$, $|g'(z)|\le c_{17}/n$
$(c_j>0$, $j=\overline{6,17})$. This implies that  
$c_{18}/n\le|(g(z)(f_n(z))'|\le c_{19}/n$ if $z=\delta_n'$ or $z=\delta_n''$.
It follows from the last inequality and (29) that for the same $z$ we have $c_{20}/n\le|F_n'(z)|\le c_{21}/n$
$(c_j>0$, $j=\overline{18,21})$. This, together with (32), yields estimate (31). Case 2 can be analyzed in a similar way.

Let us estimate the product $u_n(x)\overline{v_n(\xi)}$.
Let $H_0(x,\xi,\mu)$, $g_0(x,\xi)$, $W_0(\xi)$ and $"_0(x,\xi,\mu)$
be the functions given by (25-28) with $\varphi(x,\mu)$ and
$\psi(x,\mu)$ replaced by $e^{i\mu x}$ and $e^{-i\mu x}$.

Let us prove that
$$
H(x,\xi,\mu_n)-H_0(x,\xi,\mu_n)=O(n^{-1})\eqno(33)
$$
in case 1. Since $\mu_n=2\pi n+O(n^{-1})$, it follows from (22) that $\varphi(x,\mu_n)=e^{2\pi inx}+O(n^{-1})$,
$B_1(\varphi(x,\mu_n))=B_1(e^{2\pi inx})+O(1)$,        $B_2(\varphi(x,\mu_n))=B_2(e^{2\pi inx})+O(n^{-1})$,
$B_1^\star(\varphi(x,\mu_n))=B_1^\star(e^{2\pi inx})+O(1)$, $B_2^\star(\varphi(x,\mu_n))=B_2^\star(e^{2\pi inx})+O(n^{-1})$.
Similar estimates are valid for the functions $\psi(x,\mu_n)$ and $e^{-2\pi inx}$.
This, together with (28), yields $\Phi(x,\xi,\mu_n)=\Phi_0(x,\xi,2\pi n)+O(1)$. It also follows from (22)
that $\Phi_0(x,\xi,2\pi n)=O(n)$. It can easily be checked that $W_0(\xi,2\pi n)=4\pi in$ and
$W(\xi,\mu_n)=W_0(\xi,2\pi n)+O(1)$. The last four relations and formula (27) mean that
estimate (33) holds.

From (24) and (33) we obtain
$$
u_n(x)\overline{v_n(\xi)}=(-4\pi nH_0(x,\xi,2\pi n)+O(1))/\Delta'(\mu_n).
$$
The expression for $H_0(x,\xi,2\pi n)$ was computed in [7, c. 329]:
$$
-2\pi inH_0(x,\xi,2\pi n)=A_{34}(\cos2\pi n(x-\xi)-\cos2\pi n(x+\xi))
$$
$(x\ne\xi)$. It follows from the last two relations and (31) that
$$
|u_n(x)\overline{v_n(\xi)}|\le C.\eqno(34)
$$
By the same argument, we obtain estimate (34) in case 2. It follows from (34) [8] that the system $\{u_n(x)\}$
forms a Riesz basis in $L_2(0,1)$. Theorem 1 is proved.

It was shown in [6] that any boundary conditions (15) satisfying the requirements
of Theorem 1 are equivalent to the boundary conditions specified by the matrix 
$$
A=\left(
\begin{array}{cccc}
1&-1&0&b_0\\
0&0&1&-1
\end{array}
\right)\quad or\quad A=\left(
\begin{array}{cccc}
1&1&0&b_0\\
0&0&1&1
\end{array}
\right);
$$
in both cases, $b_0\ne0$.

{\bf Theorem 2.} {\it If
$$
A_{14}\ne A_{23},\eqno(35)
$$ then the system
of root functions $\{u_n(x)\}$ of problems (19) is a Riesz basis in $L_2(0,1)$
if and only if all but finitely many eigenvalues $\lambda_n$ 
are multiple (in other words, they are asymptotically multiple
).}

Proof. Suppose, the eigenvalues $\lambda_n$ are asymptotically multiple. It is known [9] that
the two-dimensional subspaces corresponding to the
pairwise close eigenvalues form a basis in $L_2(0,1)$, which is equivalent to an
orthogonal basis. Choosing in each of these subspases corresponding to the multiple
eigenvalues an orthonormal basis, we obtain  [10, p. 414] that
the system of root functions of problem (19), which is the union of all
orthogonal bases of mentioned subspases, is a Riesz basis in $L_2(0,1)$.

Suppose, the spectrum is not asymptotically multiple. Then there exists a subsequence of numbers such that for any number $n$ from
this subsequence
$\mu_n'\ne\mu_n''$. Let $\tilde u_n(x)$ be the eigenfunction
corresponding to an eigenvalue $\lambda_n'$ from this subsequence, and let $\tilde v_n(x)$ be the function
in the biorthogonal system corresponding to $\tilde u_n(x)$.
Let us estimate the product $\tilde u_n(x)\overline{\tilde v_n(\xi)}$.

We consider the determinant $\Delta_0(\mu)$. It follows from [6] and (21) that
$$
\Delta_0(\mu)=-2i(A_{14}+A_{23})\mu(1\mp\cos\mu)+2iA_{34}\sin\mu,
$$
where the upper sign is chosen in case 1, and the lower sign is chosen in case 2. Differentiating,
we obtain
$$
\Delta_0'(\mu)=-2i(A_{14}+A_{23})(1\mp\cos\mu\pm\mu\sin\mu)+2iA_{34}\cos\mu.
$$
It follows from the last equality and asymptotic formulas (23) and (23$'$) that 
$|\Delta_0'(\mu_n')|\le c_1\sqrt{n}$. This, together with (20), yields
$$
|\Delta'(\mu_n')|\le c_2\sqrt{n}.\eqno(36)
$$
Let $H_0(x,\xi,\mu)$, $g_0(x,\xi)$, $W_0(\xi)$ and $"_0(x,\xi,\mu)$
be the functions given by (25-28) with $\varphi(x,\mu)$ and
$\psi(x,\mu)$ replaced by $e^{i\mu x}$ and $e^{-i\mu x}$, respectively.

Let us prove that
$$
H(x,\xi,\mu_n')-H_0(x,\xi,2\pi n)=O(n^{-1/2}).\eqno(37)
$$ in case 1.
Since $\mu_n'=2\pi n+O(n^{-1/2})$, it follows from (22), that $\varphi(x,\mu_n')=e^{2\pi inx}+O(n^{-1/2})$,
$B_1(\varphi(x,\mu_n'))=B_1(e^{2\pi inx})+O(n^{1/2})$, $B_2(\varphi(x,\mu_n'))=B_2(e^{2\pi inx})+O(n^{-1/2})$,
$B_1^\star(\varphi(x,\mu_n'))=B_1^\star(e^{2\pi inx})+O(n^{1/2})$, $B_2^\star(\varphi(x,\mu_n'))=B_2^\star(e^{2\pi inx})+O(n^{-1/2})$.
Similar estimates are valid for the functions $\psi(x,\mu_n')$ and $e^{-2\pi inx}$.
This, together with (28), yields $\Phi(x,\xi,\mu_n')=\Phi_0(x,\xi,2\pi n)+O(n^{1/2})$. It also follows from (22)
that $\Phi_0(x,\xi,2\pi n)=O(n)$. It can easily be checked that $W_0(\xi,2\pi n)=4\pi in$ and
$W(\xi,\mu_n)=W_0(\xi,2\pi n)+O(n^{1/2})$. The last four relations and formula (30) mean
that estimate (37) holds.

From (24) and (37) we obtain
$$
\tilde u_n(x)\overline{\tilde v_n(\xi)}=(-4\pi nH_0(x,\xi,2\pi n)+O(n^{1/2}))/\Delta'(\mu_n').
$$
The expression for $H_0(x,\xi,2\pi n)$ was computed in [7, c. 329]:
$$\begin{array}{c}
-2\pi inH_0(x,\xi,2\pi n)=A_{34}(\cos2\pi n(x-\xi)-\cos2\pi n(x+\xi))+\\
+2\pi n[(A_{14}+A_{23}+2A_{24})\sin2\pi n(x-\xi)-(A_{14}-A_{23})\sin2\pi n(x+\xi)]
\end{array}
$$
for $x<\xi$,
$$\begin{array}{c}
-2\pi inH_0(x,\xi,2\pi n)=A_{34}(\cos2\pi n(x-\xi)-\cos2\pi n(x+\xi))+\\
+2\pi n[(A_{14}+A_{23}+2A_{13})\sin2\pi n(x-\xi)-(A_{14}-A_{23})\sin2\pi n(x+\xi)]
\end{array}
$$
for $x>\xi$.
It follows from the last three equalities, (36) and (35) that
$$
||\tilde u_n||_{L_2(0,1)}||\tilde v_n||_{L_2(0,1)}\ge C\sqrt{n},
$$
where $C>0$, and, hence, the root function system $\{u_n(x)\}$ of problem (19)
is not a basis in $L_2(0,1)$. Case 2 can be treated in a similar way.

Thus, we have established that conditions (21) and (35) reduce the question about the basis property
for the system of eigenfunctions and associated functions to the asymptotic multiplicity of the spectrum.
The presence of this property depends essentially on the particular form of the 
boundary conditions and the function $q(x)$. In the simplest case of $q(x)\equiv0$, the problem was
solved completely in [6]. Below, we cite some results of [6].

Suppose that the boundary conditions in problem (19) satisfy (21), (35), and the condition $A_{34}=0$.
We refer to such problems as problems of type (*).

They have asymptotically multiple spectrum, and  
any boundary conditions (15) satisfying the requirements mentioned above are equivalent
to the boundary conditions determined by the matrix 
$$
A=\left(
\begin{array}{cccc}
1&b_1&0&0\\
0&0&1&d_0
\end{array}
\right),
$$
where either $b_1=\mp1$, $d_0\ne1$, and $d_0\ne-1$; $d_0=\mp1$, $b_1\ne1$, and $b_1\ne-1$;
$$
A=\left(
\begin{array}{cccc}
1&\mp1&0&0\\
0&0&0&1
\end{array}
\right) \quad or\quad
A=\left(
\begin{array}{cccc}
0&1&0&0\\
0&0&1&\mp1
\end{array}
\right).
$$
The sign is always upper in case 1 and lower in case 2.

If conditions (21) and (35) hold but $A_{34}\ne0$, then the spectrum of problem (19) is asymptotically
simple,
and any boundary conditions (15) satisfying the above requirements
are equivalent to those specified by the matrix
$$
A=\left(
\begin{array}{cccc}
1&b_1&0&b_0\\
0&0&1&d_0
\end{array}
\right),
$$
where either $b_1=\mp1$, $d_0\ne1$, $d_0\ne-1$, and $b_0\ne0$; $d_0=\mp1$, $b_1\ne1$, $b_1\ne-1$, and $b_0\ne0$;
or
$$
A=\left(
\begin{array}{cccc}
1&\mp1&a_0&0\\
0&0&0&1
\end{array}
\right), \quad where \quad a_0\ne0,
$$
or
$$
A=\left(
\begin{array}{cccc}
0&1&0&b_0\\
0&0&1&\mp1
\end{array}
\right),\quad where \quad b_0\ne0.
$$
The sign is always upper in case 1 and lower in case 2.

Suppose that conditions (21) and (35) hold. Then, in the author's opinion,  of great 
interest is the problem of finding potentials
$q(x)\not\equiv0$ that ensure an asymptotically multiple spectrum
. In this relation, we mention the following results.

In [11, 12], it was established that, under the condition
$$
q(x)=q(1-x),\eqno(38)
$$
where $x\in[0,1]$, the spectrum of each of the problems
$$
Lu+\lambda u=0,\quad u'(0)=u'(1),\quad u(0)=bu(1);
$$
$$
Lu+\lambda u=0,\quad u'(0)=bu'(1),\quad u(0)=u(1),
$$
where $b\ne-1$, coincides with that of the periodic problem
$$
Lu+\lambda u=0,\quad u'(0)=u'(1),\quad u(0)=u(1),\eqno(39)
$$
and the spectrum of each of the problems
$$
Lu+\lambda u=0,\quad u'(0)+u'(1)=0,\quad u(0)+bu(1)=0;
$$
$$
Lu+\lambda u=0,\quad u'(0)+bu'(1)=0,\quad u(0)+u(1)=0,
$$
where $b\ne-1$, coincides with that of the antiperiodic problem
$$
Lu+\lambda u=0,\quad u'(0)+u'(1)=0,\quad u(0)+u(1)=0.\eqno(40)
$$
Therefore, under condition (38) the spectrum of a problem of type
(*) coincides with the spectrum of problem (39) or (40).

Let $q(x)\in L_2(0,1)$ be a real-valued function. We denote the eigenvalues of problem (39) by
 $\lambda_0$,
$\lambda_n^-$, and $\lambda_n^+$, where $n=2k$ and $k=1,2,\ldots$, and the eigenvalues
of problem (40) by $\lambda_n^-$ and $\lambda_n^+$, where $n=2k-1$ and $k=1,2,\ldots$;
in both cases, the eigenvalues are enumerated in nondecreasing order. Let $\gamma_n=\lambda_n^+-\lambda_n^-$ $(n=1,2,\ldots)$ be
the length of the spectral gap. In [13] estimates for $\gamma_n$ were obtained
for problems (39) and (40) with the potential
$$
q(x)=-\pi^2(4\alpha t\cos2\pi x+2\alpha^2\cos4\pi x),\eqno(41)
$$
where $\alpha,t$ are real numbers, and $\alpha\ne0$ and $t\ne0$.
In particular, in [13] it was shown that for even $n$
$$
\gamma_n=\frac{8\pi^2\alpha^n}{2^n[(n-2)!!]^2}|\cos(\frac{\pi}{2}t)|[1+O((\log n)^3/n)], 
$$
and for odd $n$
$$
\gamma_n=\frac{8\pi^2\alpha^n}{2^n[(n-2)!!]^2}\frac{2}{\pi}|\sin(\frac{\pi}{2}t)|[1+O((\log n)^3/n)]. 
$$
Since for potential (41) condition (38) holds, we see that for any problem of type (*) the parameter $t$ 
can be chosen so that its spectrum is asymptotically multiple or asymptotically simple.

If boundary conditions satisfy (21) and (35), then it follows from [4] that under 
supplementary conditions $q(x)\in W_1^1[0,1]$ and $2A_{34}^2\ne(A_{13}+A_{24})(A_{14}-A_{23})(q(1)-q(0))$
the spectrum of problem (19) is asymptotically simple and the root function system is not a basis.
For a problem of type (*) the last condition is equivalent to the condition $q(1)\ne q(0)$. It is
readily seen that for the potential determined by (41) $q(1)=q(0)$ for
any $\alpha$ and $t$, hence, in comparison with [4], Theorem 2 of the present paper widens the class of 
boundary value problems such that the corresponding root function system is not a basis.

It is known [14] that the spectrums of periodic and antiperiodic problems on the segment
$[0,1]$ for the Mathieu operator
$lu=u''-2\pi^2a\cos2\pi x$, where $a$ is a real number $(a\ne0)$, are 
simple. It follows from our reasoning that the eigenfunction system of the 
Mathieu operator with boundary conditions of type (*) is not a basis.
It is clear that this example is not covered by [4].
      
\medskip
\medskip

\centerline {\bf References}
\medskip
\medskip
[1] V.A. Il'in, On the basis property of the root function systems of nonselfadjoint
differential operators, in Selected Topics in Mathematics, Mechanics and
Their Applications,(Moscow, 1999), pp. 223-229.

[2] V.A. Il'in, Bases formed by root functions of differential operators,
in Program Systems (Moscow, 1999), pp. 36-43.

[3] V.A. Il'in, On a connection between the form of the boundary conditions and the basis property
and the property of equiconvergence with a trigonometric series of expansions in
root functions of a selfadjoint differential operator, Differ. Uravn., {\bf 30}, 1516-1529 (1994).

[4] A.S. Makin, A class of boundary value problems for the Sturm-Liouville operator, Differ. Uravn.,
{\bf 35}, 1058-1066 (1999).

[5] M.A. Naimark, Linear Differential Operators, (Moscow, 1969).

[6] P. Lang, J. Locker, Spectral theory of two-point differential operators
determined by $-D^2$, J. Math. Anal. Appl. {\bf 146}, 148-191 (1990).

[7] E. Coddington and N. Levinson, Theory of Ordinary Differential Equations,
(Moscow, 1958).

[8] V.A. Il'in, On the property of being an unconditional basis on a closed interval
of the system of eigenfunctions and associated functions of a second order
differential operator, Dokl. Akad. Nauk SSSR {\bf 273}, 1048-1053 (1983).

[9] A.A. Shkalikov, Bases formed by eigenfunctions of ordinary differential operatots with integral
boundary conditions, Vestnik MGU, Matem. i mekh., No. 6, 12-21 (1982).

[10] I. Ts. Gokhberg and M.G. Krein, An Introduction to the Theory of Linear Non-Self-Adjoint
Operators in Hilbert Spaces (Moscow, 1965).

[11] V.A. Sadovnichii, B.E. Kanguzhin, On a connection between the spectrum of a differential operator
with symmetric coefficients and boundary conditions, Dokl. Akad. Nauk SSSR {\bf 267}, 310-313 (1982).

[12] B.E. Kanguzhin, Some questions of the theory of inverse problems, Ph.D. Thesis, Moscow, MGU, 1982.

[13] P. Djakov, B. Mityagin, Asymptotics of instability zones of the Hill operator with a two
term potential, OSU Math. Res. Inst. Preprint 04-10 (07.10.04). P. 1-39.

[14] E.L. Ince, A proof of the impossibility of the coexistence of two Mathieu functions, Proc. Camb. Phil. Soc. {\bf 21}, 117-120 (1922).

\vspace{15mm}

Moscow State Academy of Instrument-Making and Informatics,
Stromynka 20, Moscow, 107996, Russia

\vspace{10mm}

E-mail address: alexmakin@yandex.ru

}
\end{document}